# PSEUDO-MAXIMUM LIKELIHOOD ESTIMATION OF ARCH(∞) MODELS

BY PETER M. ROBINSON[1] AND PAOLO ZAFFARONI

*London School of Economics and Imperial College London*

Strong consistency and asymptotic normality of the Gaussian pseudo-maximum likelihood estimate of the parameters in a wide class of ARCH(∞) processes are established. The conditions are shown to hold in case of exponential and hyperbolic decay in the ARCH weights, though in the latter case a faster decay rate is required for the central limit theorem than for the law of large numbers. Particular parameterizations are discussed.

**1. Introduction.** ARCH(∞) processes comprise a wide class of models for conditional heteroscedasticity in time series. Consider, for $t \in \mathbb{Z} = \{0, \pm 1, \ldots\}$, the equations

$$(1) \qquad x_t = \sigma_t \varepsilon_t,$$

$$(2) \qquad \sigma_t^2 = \omega_0 + \sum_{j=1}^{\infty} \psi_{0j} x_{t-j}^2,$$

where

$$(3) \qquad \omega_0 > 0, \qquad \psi_{0j} > 0 \qquad (j \geq 1), \qquad \sum_{j=1}^{\infty} \psi_{0j} < \infty,$$

and $\{\varepsilon_t\}$ is a sequence of independent identically distributed (i.i.d.) unobservable real-valued random variables. We shall assume that a strictly stationary solution $x_t$ to (1) and (2) exists almost surely (a.s.) under (3), and call it an ARCH(∞) process. We consider a parametric version, in which we

Received October 2003; revised January 2005.

[1]Supported by a Leverhulme Trust Personal Research Professorship and ESRC Grants R000238212 and R000239936.

know functions $\psi_j(\zeta)$ of the $r \times 1$ vector $\zeta$, for $r < \infty$, such that, for some unknown $\zeta_0$,

$$\psi_j(\zeta_0) = \psi_{0j}, \qquad j \geq 1. \tag{4}$$

Also, $\omega_0$ is unknown and $x_t$ is unobservable but we observe

$$y_t = \mu_0 + x_t \tag{5}$$

for some unknown $\mu_0$.

ARCH($\infty$) processes, extending the ARCH($m$), $m < \infty$, process of Engle [11] and the GARCH($n, m$) process of Bollerslev [4], were considered by Robinson [29] as a class of parametric alternatives in testing for serial independence of $y_t$. Empirical evidence of Whistler [35] and Ding, Granger and Engle [10] has suggested the possibility of long memory autocorrelation in the squares of financial data. Taking [contrary to the first requirement in (3)] $\omega_0 = 0$, such long memory in $x_t^2$ driven by (1) and (2) was considered by Robinson [29], the $\psi_{0j}$ being the autoregressive weights of a fractionally integrated process, implying $\sum_{j=1}^{\infty} \psi_{0j} = 1$; see also Ding and Granger [9]. For such $\psi_{0j}$, and the same objective function as was employed to generate the tests of Robinson [29], Koulikov [20] established asymptotic statistical properties of estimates of $\zeta_0$. On the other hand, under our assumption $\omega_0 > 0$, Giraitis, Kokoszka and Leipus [13] found that such $\psi_{0j}$ are inconsistent with covariance stationarity of $x_t$, which holds when $\sum_{j=1}^{\infty} \psi_{0j} < 1$. Finite variance of $x_t$ implies summability of coefficients of a linear moving average in martingale differences representation of $x_t^2$; see [37]. In this paper we do not assume finite variance of $x_t$, but rather that $x_t$ has a finite fractional moment of degree less than 2. The first requirement in (3) was shown by Kazakevičius and Leipus [18] to be necessary for existence of an $x_t$ satisfying (1) and (2). The intermediate requirement in (3) is sufficient but not necessary for a.s. positivity of $\sigma_t^2$, and is imposed here to facilitate a clearer focus on the $\psi_{0j}$, which decay, possibly slowly, but never vanish.

We wish to estimate the $(r + 2) \times 1$ vector $\theta_0 = (\omega_0, \mu_0, \zeta_0')'$ on the basis of observations $y_t$, $t = 1, \ldots, T$, the prime denoting transposition. The case when $\mu_0$ is known, for example, $\mu_0 = 0$, is covered by a simplified version of our treatment. If the $y_t$ were instead unobserved regression errors, we have $\mu_0 = 0$, but would then need to replace $x_t$ by residuals in what follows; the details of this extension would be relatively straightforward. Another relatively straightforward extension would cover simultaneous estimation of the regression parameters $\omega_0$ and $\zeta_0$, after replacing $\mu_0$ by a more general parametric function; as in (1), (2) and (5), efficiency gain is afforded by simultaneous estimation.

Under stronger restrictions than $\sum_{j=1}^{\infty} \psi_{0j} < 1$, Giraitis and Robinson [14] considered discrete-frequency Whittle estimation of $\zeta_0$, based on the squared observations $y_t^2$ (with $\mu_0$ known to be zero), this being asymptotically equivalent to constrained least squares regression of $y_t^2$ on the $y_{t-s}^2$, $s > 0$, a method employed in special cases of (2) by Engle [11] and Bollerslev [4]. In these the spectral density of $y_t^2$, when it exists, has a convenient closed form. This property, along with availability of the fast Fourier transform, makes discrete-frequency Whittle estimation based on the $y_t^2$ a computationally attractive option for point estimation, even in very long financial time series. However, it has a number of disadvantages, as discussed by Giraitis and Robinson [14]: it is not only asymptotically inefficient under Gaussian $\varepsilon_t$, but never asymptotically efficient; it requires finiteness of fourth moments of $y_t$ for consistency and of eighth moments for asymptotic normality, which are sometimes considered unacceptable for financial data; its limit covariance matrix is relatively complicated to estimate; it is less well motivated in ARCH models than in stochastic volatility and nonlinear moving average models, such as those of Taylor [33], Robinson and Zaffaroni [30, 31], Harvey [15], Breidt, Crato and de Lima [5] and Zaffaroni [36], where the actual likelihood is relatively intractable, while Whittle estimation also plays a less special role in the short-memory-in-$y_t^2$ ARCH models of Giraitis and Robinson [14] than in the long-memory-in-$y_t^2$ models of the previous five references, where it entails automatic "compensation" for possible lack of square-integrability of the spectrum of $y_t^2$. Mikosch and Straumann [26] have shown that a finite fourth moment is necessary for consistency of Whittle estimates, and that convergence rates are slowed by fat tails in $\varepsilon_t$.

For Gaussian $\varepsilon_t$, a widely-used approximate maximum likelihood estimate is defined as follows. Denote by $\theta = (\omega, \mu, \zeta')'$ any admissible value of $\theta_0$ and define

$$x_t(\mu) = y_t - \mu,$$

$$\sigma_t^2(\theta) = \omega + \sum_{j=1}^{\infty} \psi_j(\zeta) x_{t-j}^2(\mu)$$

for $t \in \mathbb{Z}$, and

$$\bar{\sigma}_t^2(\theta) = \omega + \sum_{j=1}^{t-1} \psi_j(\zeta) x_{t-j}^2(\mu) \mathbb{1}(t \geq 2)$$

for $t \geq 1$, where $\mathbb{1}(\cdot)$ denotes the indicator function. Define also

$$q_t(\theta) = \frac{x_t^2(\mu)}{\sigma_t^2(\theta)} + \ln \sigma_t^2(\theta), \qquad \bar{q}_t(\theta) = \frac{x_t^2(\mu)}{\bar{\sigma}_t^2(\theta)} + \ln \bar{\sigma}_t^2(\theta), \qquad 1 \leq t \leq T,$$

$$Q_T(\theta) = T^{-1} \sum_{t=1}^{T} q_t(\theta), \qquad \bar{Q}_T(\theta) = T^{-1} \sum_{t=1}^{T} \bar{q}_t(\theta),$$

$$\tilde{\theta}_T = \arg\min_{\theta \in \Theta} Q_T(\theta), \qquad \hat{\theta}_T = \arg\min_{\theta \in \Theta} \bar{Q}_T(\theta),$$

where $\Theta$ is a prescribed compact subset of $\mathbb{R}^{r+2}$. The quantities with over-bar are introduced due to $y_t$ being unobservable for $t \leq 0$; $\tilde{\theta}_T$ is uncomputable. Because we do not assume Gaussianity in the asymptotic theory, we refer to $\hat{\theta}_T$ as a pseudo-maximum likelihood estimate (PMLE).

We establish strong consistency of $\hat{\theta}_T$ and asymptotic normality of $T^{1/2}(\hat{\theta}_T - \theta_0)$, as $T \to \infty$, for a class of $\psi_j(\zeta)$ sequences. In the case of the first property this is accomplished by first showing strong consistency of $\tilde{\theta}_T$ and then that $\hat{\theta}_T - \tilde{\theta}_T \to 0$, a.s. In the case of the second we likewise first show it for $T^{1/2}(\tilde{\theta}_T - \theta_0)$ and then show that $\hat{\theta}_T - \tilde{\theta}_T = o_p(T^{-1/2})$, but the latter property, and thus the asymptotic normality of $T^{1/2}(\hat{\theta}_T - \theta_0)$, is achieved only under a restricted set of possible $\zeta_0$ values, and this seems of practical concern in relation to some popular choices of the $\psi_j(\zeta)$. These results are presented in the following section, along with a description of regularity conditions and partial proof details. The structure of the proof is similar in several respects to earlier ones for the GARCH case of (2), especially that of Berkes, Horváth and Kokoszka [3]. Sections 3 and 4 apply the results to particular models.

## 2. Assumptions and main results. Our assumptions are as follows.

ASSUMPTION A($q$), $q \geq 2$.   The $\varepsilon_t$ are i.i.d. random variables with $E\varepsilon_0 = 0$, $E\varepsilon_0^2 = 1$, $E|\varepsilon_0|^q < \infty$ and probability density function $f(\varepsilon)$ satisfying

$$f(\varepsilon) = O(L(|\varepsilon|^{-1})|\varepsilon|^b) \qquad \text{as } \varepsilon \to 0,$$

for $b > -1$ and a function $L$ that is slowly varying at the origin.

ASSUMPTION B.   There exist $\omega_L, \omega_U, \mu_L, \mu_U$ such that $0 < \omega_L < \omega_U < \infty$, $-\infty < \mu_L < \mu_U < \infty$, and a compact set $\Upsilon \in R^r$ such that $\Theta = [\omega_L, \omega_U] \times [\mu_L, \mu_U] \times \Upsilon$.

ASSUMPTION C.   $\theta_0$ is an interior point of $\Theta$.

ASSUMPTION D.   For all $j \geq 1$,

$$\inf_{\zeta \in \Upsilon} \psi_j(\zeta) > 0; \tag{6}$$

$$\sup_{\zeta \in \Upsilon} \psi_j(\zeta) \leq K j^{-\underline{d}-1} \qquad \text{for some } \underline{d} > 0; \tag{7}$$

$$\psi_{0j} \leq K \psi_{0k} \qquad \text{for } 1 \leq k \leq j, \tag{8}$$

where $K$ throughout denotes a generic, positive constant.

ASSUMPTION E.  There exists a strictly stationary and ergodic solution $x_t$ to (1) and (2), and for some

$$\rho \in ((\underline{d}+1)^{-1}, 1), \tag{9}$$

with $\underline{d}$ as in Assumption D, we have

$$E|x_0|^{2\rho} < \infty. \tag{10}$$

ASSUMPTION F ($l$).  For all $j \geq 1$, $\psi_j(\zeta)$ has continuous $k$th derivative on $\Upsilon$ such that, with $\zeta_i$ denoting the $i$th element of $\zeta$,

$$\left| \frac{\partial^k \psi_j(\zeta)}{\partial \zeta_{i_1} \cdots \partial \zeta_{i_k}} \right| \leq K \psi_j(\zeta)^{1-\eta} \tag{11}$$

for all $\eta > 0$ and all $i_h = 1, \ldots, r$, $h = 1, \ldots, k$, $k \leq l$.

ASSUMPTION G.  For each $\zeta \in \Upsilon$, there exist integers $j_i = j_i(\zeta)$, $i = 1, \ldots, r$, such that $1 \leq j_1(\zeta) < \cdots < j_r(\zeta) < \infty$ and

$$\mathrm{rank}\{\Psi_{(j_1, \ldots, j_r)}(\zeta)\} = r,$$

where

$$\Psi_{(j_1, \ldots, j_r)}(\zeta) = \{\psi_{j_1}^{(1)}(\zeta), \ldots, \psi_{j_r}^{(1)}(\zeta)\}, \qquad \psi_j^{(1)}(\zeta) = \frac{\partial \psi_j(\zeta)}{\partial \zeta}.$$

ASSUMPTION H.  There exists

$$d_0 > \tfrac{1}{2} \tag{12}$$

such that

$$\psi_{0j} \leq K j^{-1-d_0}, \tag{13}$$

and (10) holds for

$$\rho \in (4/(2d_0 + 3), 1). \tag{14}$$

Assumption A($q$) allows some asymmetry in $\varepsilon_t$, but implies the less primitive condition (which does not even require existence of a density) employed in a similar context by Berkes, Horváth and Kokoszka [3]. Assumptions B and C are standard. Inequalities (7) and (13) together imply $d_0 \geq \underline{d}$, while (8) with (3) is milder than monotonicity but implies $\psi_{0j} = o(j^{-1})$ as $j \to \infty$. We take $\eta > 0$ in Assumption F($l$) because $\psi_j(\zeta) < 1$ for all large enough $j$, by (7). Assumption G is crucial to the proof of consistency, being used in Lemmas 9 and 10 to show that in the limit $\theta_0$ globally minimizes $Q_T(\theta)$; it also ensures nonsingularity of the matrix $H_0$ in Proposition 2 and Theorem 2 below. This and other assumptions are discussed in Sections 3 and 4 in connection with some parameterizations of interest.

We present asymptotic results for the uncomputable $\tilde{\theta}_T$ as propositions, those for $\hat{\theta}_T$ as theorems. All these, and the corollaries in Sections 3 and 4 and lemmas in Section 5, assume (1)–(5).

PROPOSITION 1. *For some $\delta > 0$, let Assumptions* A*($2 + \delta$),* B, C, D, E, F*(1) and* G *hold. Then*

$$\tilde{\theta}_T \to \theta_0 \qquad a.s. \ as \ T \to \infty.$$

PROOF. The proof follows as in, for example, [17], Theorem 6, from uniform a.s. convergence over $\Theta$ of $Q_T(\theta)$ to $Q(\theta) = Eq_0(\theta)$ established in Lemma 7, the fact that $Q_T(\hat{\theta}_T) \leq Q_T(\theta)$, and Lemma 10. $\square$

THEOREM 1. *For some $\delta > 0$, let Assumptions* A($2 + \delta$), B, C, D, E, F*(1) and* G *hold. Then*

(15) $$\hat{\theta}_T \to \theta_0 \qquad a.s. \ as \ T \to \infty.$$

PROOF. From Lemmas 7 and 8, $\bar{Q}_T(\theta)$ converges uniformly to $Q(\theta)$ a.s., whence the proof is as indicated for Proposition 1. $\square$

Denote by $\kappa_j$ the $j$th cumulant of $\varepsilon_t$ and introduce

$$G_0 = (2 + \kappa_4)M - 2\kappa_3(N + N') + P, \qquad H_0 = M + \tfrac{1}{2}P,$$

where

$$M = E(\tau_0 \tau_0'), \qquad N = E(\sigma_0^{-1}\tau_0)e_2', \qquad P = E(\sigma_0^{-2})e_2 e_2',$$

for $\tau_0 = \tau_0(\theta_0)$, $\tau_t(\theta) = (\partial/\partial\theta)\log\sigma_t^2(\theta)$, and $e_2$ the second column of the $(r + 2) \times (r + 2)$ identity matrix. In case $\mu_0$ is known (e.g., to be zero), we omit the second row and column from $M$, and have instead $G_0 = (2 + \kappa_4)M$, $H_0 = M$. In case $\varepsilon_t$ is Gaussian, $\kappa_3 = \kappa_4 = 0$, so $G_0 = 2H_0 = 2M + P$.

PROPOSITION 2. *Let Assumptions* A(4), B, C, D, E, F(3) *and* G *hold. Then*

$$T^{1/2}(\tilde{\theta}_T - \theta_0) \xrightarrow{d} N(0, H_0^{-1}G_0 H_0^{-1}) \qquad as \ T \to \infty.$$

PROOF. Write

$$Q_T^{(1)}(\theta) = \frac{\partial Q_T(\theta)}{\partial\theta} = T^{-1}\sum_{t=1}^{T} u_t(\theta),$$

where

$$u_t(\theta) = \tau_t(\theta)(1 - \chi_t^2(\theta)) + \sigma_t^{-2}(\theta)\nu_t(\theta),$$

with

$$\chi_t(\theta) = \frac{x_t^2(\mu)}{\sigma_t^2(\theta)}, \qquad \nu_t(\theta) = \frac{\partial x_t^2(\mu)}{\partial \theta} = -2x_t(\mu)e_2.$$

By the mean value theorem,

$$(16) \qquad 0 = Q_T^{(1)}(\tilde{\theta}_T) = Q_T^{(1)}(\theta_0) + \tilde{H}_T(\tilde{\theta}_T - \theta_0),$$

where $\tilde{H}_T$ has as its $i$th row the $i$th row of $H_T(\theta) = T^{-1}\sum_{t=1}^{T} h_t(\theta)$ evaluated at $\theta = \tilde{\theta}_T^{(i)}$, where $h_t(\theta) = (\partial^2/\partial\theta\,\partial\theta')Q_T(\theta)$, $\|\tilde{\theta}_T^{(i)} - \theta_0\| \le \|\tilde{\theta}_T - \theta_0\|$, where we define $\|A\| = \{\mathrm{tr}(A'A)\}^{1/2}$ for any real matrix $A$. Now $u_t(\theta_0) = \tau_t(\theta_0)(1 - \varepsilon_t^2) - 2e_2\varepsilon_t/\sigma_t$ is, by Lemmas 2, 3 and 7, a stationary ergodic martingale difference vector with finite variance, so from Brown [6] and the Cramér–Wold device, $T^{1/2}Q_T^{(1)}(\theta_0) \to_d N(0, G_0)$ as $T \to \infty$. Finally, by Lemma 7 and Theorem 1, $\tilde{H}_T \to_p H_0$, whence the proof is completed in standard fashion. □

Define

$$\bar{u}_t(\theta) = \frac{\partial \bar{q}_t(\theta)}{\partial \theta}, \qquad \bar{g}_t(\theta) = \bar{u}_t(\theta)\bar{u}_t'(\theta), \qquad \bar{h}_t(\theta) = \frac{\partial^2 \bar{q}_t(\theta)}{\partial\theta\,\partial\theta'},$$

$$\bar{G}_T(\theta) = T^{-1}\sum_{t=1}^{T} \bar{g}_t(\theta), \qquad \bar{H}_T(\theta) = T^{-1}\sum_{t=1}^{T} \bar{h}_t(\theta).$$

THEOREM 2. *Let Assumptions* A(4), B, C, D, E, F(3), G *and* H *hold. Then*

$$(17) \qquad T^{1/2}(\hat{\theta}_T - \theta_0) \xrightarrow{d} N(0, H_0^{-1}G_0H_0^{-1}) \qquad as\ T \to \infty,$$

*and* $H_0^{-1}G_0H_0^{-1}$ *is strongly consistently estimated by* $\bar{H}_T^{-1}(\hat{\theta}_T)\bar{G}_T(\hat{\theta}_T)\bar{H}_T^{-1}(\hat{\theta}_T)$.

PROOF. We have

$$0 = \bar{Q}_T^{(1)}(\hat{\theta}_T) = \bar{Q}_T^{(1)}(\theta_0) + \hat{H}_T(\hat{\theta}_T - \theta_0),$$

where $\bar{Q}_T^{(1)}(\theta) = (\partial/\partial\theta)Q_T(\theta)$ and $\hat{H}_T$ has as its $i$th row the $i$th row of $\bar{H}_T(\theta)$ evaluated at $\theta = \hat{\theta}_T^{(i)}$, for $\|\hat{\theta}_T^{(i)} - \theta_0\| \le \|\hat{\theta}_T - \theta_0\|$. Thus, from (16),

$$\hat{\theta}_T - \tilde{\theta}_T = (\tilde{H}_T^{-1} - \hat{H}_T^{-1})\bar{Q}_T^{(1)}(\theta_0) - \tilde{H}_T^{-1}\{\bar{Q}_T^{(1)}(\theta_0) - Q_T^{(1)}(\theta_0)\},$$

where the inverses exist a.s. for all sufficiently large $T$ by Lemma 9. In view of Proposition 2 and Lemma 8, (17) follows on showing that

$$\bar{Q}_T^{(1)}(\theta_0) - Q_T^{(1)}(\theta_0) = o_p(T^{-1/2}).$$

The left-hand side can be written $(B_{1T} + B_{2T} + B_{3T})/T$, where

$$B_{1T} = \sum_{t=1}^{T} \varepsilon_t^2 b_{1t}, \qquad B_{2T} = -\sum_{t=1}^{T} (\varepsilon_t^2 - 1) b_{2t}, \qquad B_{3T} = -2e_2 \sum_{t=1}^{T} \varepsilon_t b_{3t},$$

with

$$b_{1t} = -\frac{\bar{\sigma}_t^{2(1)}(\sigma_t^2 - \bar{\sigma}_t^2)}{\bar{\sigma}_t^4}, \qquad b_{2t} = \frac{\sigma_t^{2(1)}}{\sigma_t^2} - \frac{\bar{\sigma}_t^{2(1)}}{\bar{\sigma}_t^2}, \qquad b_{3t} = \frac{\sigma_t^2 - \bar{\sigma}_t^2}{\bar{\sigma}_t^2 \sigma_t},$$

for $\bar{\sigma}_t^2 = \bar{\sigma}_t^2(\theta_0)$, $\sigma_t^{2(1)} = \sigma_t^{2(1)}(\theta_0)$, $\bar{\sigma}_t^{2(1)} = \bar{\sigma}_t^{2(1)}(\theta_0)$, with $\sigma_t^{2(1)}(\theta) = (\partial/\partial\theta)\sigma_t^2(\theta)$, $\bar{\sigma}_t^{2(1)}(\theta) = (\partial/\partial\theta)\bar{\sigma}_t^2(\theta)$. We show that $B_{iT} = o_p(T^{1/2})$, $i = 1, 2, 3$. For the remainder of this proof, we drop the zero subscript in $\psi_{0j}$.

Consider first $B_{1T}$. We have

$$\bar{\sigma}_t^{2(1)} = \left(1, -2\sum_{j=1}^{t-1} \psi_j x_{t-j}, \sum_{j=1}^{t-1} \psi_j^{(1)} x_{t-j}^2\right)', \tag{18}$$

where $\psi_j^{(1)} = \psi_j^{(1)}(\zeta_0)$. From Assumption F(1),

$$\|\bar{\sigma}_t^{2(1)}\| \leq 1 + 2\sum_{j=1}^{t-1} \psi_j |x_{t-j}| + K\sum_{j=1}^{t-1} \psi_j^{1-\eta} x_{t-j}^2,$$

for all $\eta > 0$. Now

$$\sum_{j=1}^{t-1} \psi_j |x_{t-j}| \leq \left(\sum_{j=1}^{t-1} \psi_j x_{t-j}^2\right)^{1/2} \left(\sum_{j=1}^{\infty} \psi_j\right)^{1/2} \leq K\bar{\sigma}_t,$$

so since $\bar{\sigma}_t \geq \omega_L > 0$,

$$\bar{\sigma}_t^{-2} \sum_{j=1}^{t-1} \psi_j |x_{t-j}| \leq K\bar{\sigma}_t^{-1} < \infty.$$

From (8),

$$\sum_{j=1}^{t-1} \psi_j^{1-\eta} x_{t-j}^2 \leq K\psi_t^{-\eta} \bar{\sigma}_t^2.$$

It follows that

$$\|\bar{\sigma}_t^{2(1)}\|/\bar{\sigma}_t^2 \leq K\psi_t^{-\eta}. \tag{19}$$

On the other hand, by the $c_r$-inequality ([23], page 157) and (10),

$$E(\sigma_t^2 - \bar{\sigma}_t^2)^\rho \leq K\sum_{j=t}^{\infty} \psi_j^\rho E|x_{t-j}|^{2\rho} \leq K\sum_{j=t}^{\infty} \psi_j^\rho. \tag{20}$$

Thus, by (8) and (14),

$$(21) \qquad E\|b_{1t}\|^{\rho} \leq K\psi_t^{-\eta\rho} \sum_{j=t}^{\infty} \psi_j^{\rho} \leq K \sum_{j=t}^{\infty} \psi_j^{\rho(1-\eta)} \leq Kt^{1-\rho(d_0+1)(1-\eta)},$$

choosing $\eta < 1 - 1/\{\rho(d_0+1)\}$, which (14) enables. Applying the $c_r$-inequality again,

$$E\|B_{1T}\|^{\rho} \leq K \sum_{t=1}^{T} E|\varepsilon_0|^{2\rho} E\|b_{1t}\|^{\rho}.$$

Applying (21), this is $O(1)$ when $\rho > 2/(d_0+1)$, while when $\rho \leq 2/(d_0+1)$, we may choose $\eta$ so small to bound it by

$$KT^{2-\rho(d_0+1)(1-\eta)} \leq KT^{\rho/2 - \{1+2(d_0+1)(1-\eta)\}[\rho/2 - 2/\{1+2(d_0+1)(1-\eta)\}]} = o(T^{\rho/2}),$$

using (12) [which requires (13)] and arbitrariness of $\eta$. Thus, $B_{1T} = o_p(T^{1/2})$ by Markov's inequality.

Consider $B_{2T}$. By independence of $\varepsilon_t$ and $b_{2t}$, by the $c_r$-inequality when $\rho \leq \frac{1}{2}$, and by the inequality of von Bahr and Esseen [34] and the fact that the $\varepsilon_t^2$ are i.i.d. with mean 1 when $\rho > \frac{1}{2}$,

$$E\|B_{2T}\|^{2\rho} \leq K \sum_{t=1}^{T} (E|\varepsilon_0|^{4\rho} + 1) E\|b_{2t}\|^{2\rho} \leq K \sum_{t=1}^{T} (E\|b_{4t}\|^{2\rho} + E\|b_{5t}\|^{2\rho}),$$

where

$$b_{4t} = \frac{\sigma_t^{2(1)} - \bar{\sigma}_t^{2(1)}}{\sigma_t^2}, \qquad b_{5t} = \frac{\bar{\sigma}_t^{2(1)}(\sigma_t^2 - \bar{\sigma}_t^2)}{\bar{\sigma}_t^2 \sigma_t^2}.$$

Thus, from Assumptions F(1) and H,

$$\|b_{4t}\| \leq \left( 2 \sum_{j=t}^{\infty} \psi_j |x_{t-j}| + \sum_{j=t}^{\infty} \|\psi_j^{(1)}\| x_{t-j}^2 \right) \Big/ \sigma_t^2$$

$$\leq \sigma_t^{-2} \left[ 2 \left\{ \sum_{j=t}^{\infty} \psi_j \right\}^{1/2} + \left\{ \sum_{j=t}^{\infty} (\|\psi_j^{(1)}\|^2/\psi_j) x_{t-j}^2 \right\}^{1/2} \right] \left\{ \sum_{j=t}^{\infty} \psi_j x_{t-j}^2 \right\}^{1/2}$$

$$\leq K \left\{ \left( \sum_{j=t}^{\infty} j^{-d_0-1} \right)^{1/2} + \left( \sum_{j=t}^{\infty} \psi_j^{1-2\eta} x_{t-j}^2 \right)^{1/2} \right\}$$

$$\leq K \left[ t^{-d_0/2} + \left\{ \sum_{j=t}^{\infty} j^{-(d_0+1)(1-2\eta)} x_{t-j}^2 \right\}^{1/2} \right],$$

so

$$E\|b_{4t}\|^{2\rho} \leq Kt^{-\rho d_0} + K \sum_{j=t}^{\infty} j^{-(d_0+1)\rho(1-2\eta)} \leq Kt^{1-(d_0+1)\rho(1-2\eta)}$$

for sufficiently small $\eta$. Thus, $\sum_{t=1}^{T} E\|b_{4t}\|^{2\rho}$ is $O(1)$ for $\rho > 2/(d_0+1)$, while for $\rho \leq 2/(d_0+1)$, it is bounded by

$$KT^{2-(d_0+1)\rho(1-2\eta)} \leq KT^{\rho-(d_0+2)\{\rho-2/(d_0+2)\}+2(d_0+1)\rho\eta} = o(T^\rho)$$

from (14) and arbitrariness of $\eta$. Also, $\|b_{5t}\| \leq K\|\bar{\sigma}_t^{2(1)}/\bar{\sigma}_t^2\|(\sigma_t^2-\bar{\sigma}_t^2)^{1/2}$, so from (19) and (20) we have $E\|b_{5t}\|^{2\rho} \leq Kt^{1-(d_0+1)\rho(1-2\eta)}$, and proceeding as before,

$$\sum_{t=1}^{T} E\|b_{5t}\|^{2\rho} = o(T^\rho),$$

and thence, $B_{2T} = o_p(T^{1/2})$.

Next,

$$E\|B_{3T}\|^{2\rho} \leq KE\left|\sum_{t=1}^{T} \varepsilon_t b_{3t}\right|^{2\rho} \leq K\sum_{t=1}^{T} E|\varepsilon_0|^{2\rho} E b_{3t}^{2\rho},$$

applying the $c_r$-inequality when $\rho \leq \frac{1}{2}$ and von Bahr and Esseen [34] when $\rho > \frac{1}{2}$. Now $b_{3t} \leq (\sigma_t^2 - \bar{\sigma}_t^2)^{1/2}\bar{\sigma}_t^{-2}$, so from (20),

$$E\|B_{3T}\|^{2\rho} \leq K\sum_{t=1}^{T} \sum_{j=t}^{\infty} \psi_j^\rho$$

$$\leq K\{\mathbb{1}(\rho > 2/(d_0+1)) + (\ln T)\mathbb{1}(\rho = 2/(d_0+1))$$

$$+ T^{2-\rho(d_0+1)}\mathbb{1}(\rho < 2/(d_0+1))\}$$

$$= o(T^\rho),$$

much as before. Thence, $B_{3T} = o_p(T^{1/2})$.

It remains to consider the last statement of the theorem, which follows on standard application of Propositions 1 and 2, Theorem 1 and Lemmas 7 and 8. □

In earlier versions of this paper we checked the conditions in the case of GARCH$(n,m)$ models in which the $\psi_j(\zeta)$ decay exponentially and we allow the possibility that the GARCH coefficients lie in a subspace of dimension less than $m+n$; the details are available from the authors on request. However, the literature on asymptotic theory for estimates of GARCH models is now extensive, recent references including [3, 7, 12, 16, 22, 32], along with investigations of the properties of the models themselves; see recently [2, 18, 25]. We focus instead on alternative models which have received less attention, and for which our theoretical framework is primarily intended.

We introduce the generating function

$$(22) \qquad \psi(z;\zeta) = \sum_{j=1}^{\infty} \psi_j(\zeta)z^j, \qquad |z| \leq 1.$$

**3. Fractional GARCH models.** A slowly decaying class of ARCH($\infty$) weights was considered by Robinson [29], Ding and Granger [9] and Koulikov [20], generated by

$$(23) \qquad \psi(z; \zeta) = 1 - (1-z)^\zeta, \qquad 0 < \zeta < 1,$$

where $r = 1$ and formally

$$(24) \qquad (1-z)^d = \sum_{j=0}^{\infty} \frac{\Gamma(j-d)}{\Gamma(-d)\Gamma(j+1)} z^j, \qquad |z| \le 1, \ d > 0.$$

In these references $\omega_0 = 0$ was assumed in (2), but we assume $\omega_0 > 0$ and generalize (23) as follows. Introduce the functions $a_j = a_j(\zeta)$, $b_j = b_j(\zeta)$ and, for $m \ge 1$, $n \ge 0$, $n + m \ge r$,

$$(25) \qquad a(z; \zeta) = \sum_{j=1}^{m} a_j z^j, \qquad b(z; \zeta) = 1 - \sum_{j=1}^{n} b_j z^j \mathbb{1}(n \ge 1);$$

and for all $\zeta \in \Upsilon$,

$$(26) \qquad a_j > 0, \qquad j = 1, \ldots, m; \qquad b_j > 0, \qquad j = 1, \ldots, n;$$

$$(27) \qquad b(z; \zeta) \ne 0, \qquad |z| \le 1;$$

$$(28) \qquad a(z; \zeta) \text{ and } b(z; \zeta) \text{ have no common zeros in } z.$$

Now take $\psi(z; \zeta)$ (22) to be given by

$$(29) \qquad \psi(z; \zeta) = \frac{a(z; \zeta)\{1 - (1-z)^d\}}{z b(z; \zeta)},$$

with $d = d(\zeta)$ satisfying

$$(30) \qquad d \in (0, 1).$$

We call $x_t$ based on (29) a fractional GARCH, FGARCH($n, d_0, m$) process, for $d_0 = d(\zeta_0)$.

COROLLARY 1. *Let $\psi(z; \zeta)$ be given by (29) and (25) with $m \ge 1$, $n \ge 0$, and let $d$ and the $a_j, b_j$ be continuously differentiable. For some $\delta > 0$, let Assumptions* A($2+\delta$), B, C *and* E *hold, with all $\zeta \in \Upsilon$ satisfying (26)–(28), (30) and*

$$\text{rank}\left\{ \frac{\partial}{\partial \zeta}(a_1, \ldots, a_m, b_1, \ldots, b_n, d) \right\} = r.$$

*Then (15) is true. Let also $d$ and the $a_j, b_j$ be thrice continuously differentiable and $d_0 > \frac{1}{2}$. Then (17) is true.*

PROOF. Denoting by $c_j$ $(j \geq 1)$ and $d_j$ $(j \geq 0)$ the coefficients of $z^j$ in the expansions of $a(z; \zeta)/b(z; \zeta)$, $z^{-1}\{1 - (1 - z)^d\}$, respectively, we have $\psi_j(\zeta) = \sum_{k=0}^{j-1} c_{j-k} d_k$, $j \geq 1$. From [3], the $c_j$ are bounded above and below by positive, exponentially decaying sequences when $n \geq 1$, and are all nonnegative when $n = 0$. Since the $d_j$ are all positive, it follows that (6) holds. Also, Stirling's approximation indicates that $j^{-d-1}/K \leq d_j \leq Kj^{-d-1}$, so the $\psi_j(\zeta)$ satisfy the same inequalities. Compactness of $\Upsilon$, smoothness of $d$, and (30), imply $d(\zeta) \geq \underline{d}$, to check (7). The above argument indicates that $\psi_{0j} \leq Kj^{-d_0-1} \leq Kk^{-d_0-1} \leq K\psi_{0k}$ for $j > k \geq 1$, so (8) holds, and thus Assumption D. With regard to (11), note that $(\partial/\partial d)\psi(z; \zeta) = -\{a(z; \zeta)/b(z; \zeta)\}z^{-1}(1-z)^d \ln(1-z)$, where the coefficient of $z^j$ in $-z^{-1}(1-z)^d \ln(1-z)$ is $\sum_{k=1}^{j} k^{-1} d_{j-k} \leq K(\ln j)j^{-d-1} \leq Kj^{-(d+1)(1-\eta)} \leq K\psi_j^{1-\eta}(\zeta)$ for any $\eta > 0$. Derivatives with respect to the $a_j, b_j$ are dominated, and higher derivatives can be dealt with similarly, to complete the checking of Assumption F($l$). To check Assumption G, suppress reference to $\zeta$ in $a$, $b$, $\psi$ and

$$\phi(z) = b(z)^{-1}\{1 - (1 - z)^d\}, \qquad \gamma(z) = b(z)^{-1}a(z),$$

and note that

$$\frac{\partial \psi(z)}{\partial a_j} = z^{j-1}\phi(z), \qquad j = 1, \ldots, m,$$

$$\frac{\partial \psi(z)}{\partial b_j} = z^{j-1}\gamma(z)\phi(z), \qquad j = 1, \ldots, n,$$

$$\frac{\partial \psi(z)}{\partial d} = -\frac{\gamma(z)}{z}(1 - z)^d \log(1 - z).$$

Choose $j_i(\zeta) = i$ for $i = 1, \ldots, m + n$, $\zeta \in \Upsilon$, leaving $j_{m+n+1}(\zeta)$ to be determined subsequently. Fix $\zeta$ and write $U = \Psi_{(j_i, \ldots, j_r)}(\zeta)$, partitioning it in the ratio $m + n : 1$ and calling its $(i, j)$th submatrix $U_{ij}$. We first show that the $(m + n) \times (m + n)$ matrix $U_{11}$ is nonsingular. Write $R$ for the $n \times (m + n)$ matrix with $(i, j)$th element $\gamma_{j-i}$, and $S$ for the $(m + n) \times (m + n)$ matrix with $(i, j)$th element $\phi_{j-i+1}$, where $\phi_j = \gamma_j = 0$ for $j \leq 0$, and for $j > 0$, $\phi_j$ and $\gamma_j$ are respectively given by

$$\phi(z) = \sum_{j=1}^{\infty} \phi_j z^j, \qquad \gamma(z) = \sum_{j=1}^{\infty} \gamma_j z^j,$$

these series converging absolutely for $|z| \leq 1$ in view of (30). Noting that $\psi_j^{(1)}$ is given by $(\partial/\partial \zeta)\psi(z) = \sum_{j=1}^{\infty} \psi_j^{(1)} z^j$, we find that the first $m$ rows of $U_{11}$ can be written $(I_m, O)S$, where $I_m$ is the $m$-rowed identity matrix, $O$ is the $m \times n$ matrix of zeroes and, when $n \geq 1$ the last $n$ rows of $U_{11}$ can be written $RS$. Now $S$ is upper-triangular with nonzero diagonal elements.

Thus, for $n = 0$, $U_{11} = S$ is nonsingular. For $n \geq 1$, $U_{11}$ is nonsingular if and only if the $n \times n$ matrix $R_2$ having $(i,j)$th element $\gamma_{m+j-i}$ and consisting of the last $n$ columns of $R$ is nonsingular. This is not so if and only if the $\gamma_j$, $j = m, \ldots, m+n-1$, are generated by a homogeneous linear difference equation of degree $n-1$, that is, if there exist scalars $\lambda_0, \lambda_1, \ldots, \lambda_{n-1}$, not all zero, such that

$$\lambda_0 \gamma_j - \sum_{i=1}^{n-1} \lambda_i \gamma_{j-i} = 0, \qquad j = m, \ldots, m+n-1.$$

But it follows from (25) and (27) that they are generated by the linear difference equation

$$\gamma_j - \sum_{i=1}^{n-1} b_i \gamma_{j-i} = \pi_j, \qquad j = m, \ldots, m+n-1,$$

where $\pi_m = a_m + b_n \gamma_{m-n}$, $\pi_j = b_n \gamma_{j-n}$ for $j = m+1, \ldots, m+n-1$. Since $b_n \neq 0$, the $\pi_j$ are all zero if and only if $\gamma_{m-n} = -a_m/b_n$ and $\gamma_j = 0$ for $j = m+1-n, \ldots, m-1$. But this implies $\gamma_m = 0$ also, and thence, $\gamma_j = 0$, all $j \geq m-n+1$. For $m \leq n$, this is inconsistent with the requirement $a_j > 0$, $j = 1, \ldots, m$, and for $m > n$, it implies $a$ has a factor $b$, which is inconsistent with (28). Thus, $U_{11}$ is nonsingular when $n \geq 1$. Nonsingularity of $U$ follows if $U_{22} \neq U_{21} U_{11}^{-1} U_{12}$. For large enough $j_{m+n+1} = j_{m+n+1}(\zeta)$, this must be true because $U_{22}$ decays like $(\ln j_{m+n+1}) j_{m+n+1}^{-d-1}$, whereas the elements of $U_{12}$ are $O(\beta^{j_{m+n+1}})$ for some $\beta \in (0,1)$. Thus Assumption G is true, and thence (15). Clearly (13) is true, so under the additional conditions so is Assumption H, and thence (17). $\square$

For $m = 1$, $n = 0$, (29) reduces to (23) when $a_1 = 1$, while when $a_1 \in (0,1)$, it gives model (4.24) of Ding and Granger [9]. The important difference between these two cases is that the covariance stationarity condition $\psi(1; \zeta_0) < 1$ is satisfied in the second but not in the first. In general with (29), as with the GARCH model, $x_t$ is covariance stationary when $a(1; \zeta_0) < b(1; \zeta_0)$ but not otherwise. We compare (29) with

$$(31) \qquad \psi(z; \zeta) = 1 - \frac{\{1 - a(z; \zeta)\}}{b(z; \zeta)} (1-z)^d,$$

with $d$ again satisfying (30) and $a$ and $b$ again given as in (25), though we now allow $m = 0$, meaning $a(z; \zeta) \equiv 0$. Thus, with $m = n = 0$, (31) reduces to (23). ARCH($\infty$) models with $\psi$ given by (31) were proposed by Baillie, Bollerslev and Mikkelsen [1] and called FIGARCH($n, d_0, m$). In general, though (31) also gives hyperbolically decaying $\psi_{0j}$, it differs in some notable respects. Application of (26)–(28) again ensures positivity of $\psi_j(\zeta)$ in case of FGARCH and facilitates the above proof, but sufficient conditions in FIGARCH are less apparent in general, though Baillie, Bollerslev and Mikkelsen [1] indicated that they can be obtained. Also, unlike FGARCH, FIGARCH $x_t$ never has finite variance.

The requirement $d_0 > \frac{1}{2}$ for the central limit theorem in Corollary 1 would also be imposed in a corresponding result for FIGARCH. This is automatically satisfied in GARCH models but if only $d_0 \in (0, \frac{1}{2}]$ in (13) is possible in the general setting of Section 3, it appears that the asymptotic bias in $\hat{\theta}_T$ is of order at least $T^{-1/2}$, whereas that for $\tilde{\theta}_T$ is always $o(T^{-1/2})$. Assumption H copes with the replacement of $\sigma_t^2(\theta)$ by $\bar{\sigma}_t^2(\theta)$, the truncation error varying inversely with $d_0$. Inspection of the proof of Theorem 2 indicates that this bias problem is due to the term $H^{-1}B_{1T}$. The factor $\sigma_t^2 - \bar{\sigma}_t^2$ in $b_{1t}$ is nonnegative, and if $j^{-d_0-1}$ is an exact rate for $\psi_{0j}$, $\sigma_t^2 - \bar{\sigma}_t^2$ exceeds $t^{-d_0}/K$ as $t \to \infty$ with probability approaching one. So far as the factor $\bar{\sigma}_t^{2(1)}/\sigma_t^4$ in $b_{1t}$ is concerned, the second element of $\bar{\sigma}^{2(1)}$ [see (18)] has zero mean, but the first is positive, and though the $\psi_j^{(1)}$ can have elements of either sign, whenever $d_0 \le \frac{1}{2}$ it seems unlikely that the last $r$ elements of $B_{1T}$ can be $o_p(T^{1/2})$. Nor is there scope for relaxing (12) by strengthening other conditions. With regard to implications for choice of $\rho$, when $d_0 \ge 2\underline{d} + \frac{1}{2}$, (14) entails no restriction over (9).

Though results of Giraitis, Kokoszka and Leipus [13] indicate existence of a stationary solution of (1)–(3) when $\psi(1; \zeta_0) < 1$, Kazakevičius and Leipus [19] have questioned the existence of strictly stationary FIGARCH processes, and thus the relevance of Assumption E here. The same reservations can be expressed about FGARCH when $a(1; \zeta_0) \ge b(1; \zeta_0)$, and more generally about ARCH($\infty$) processes with $\psi(1; \zeta_0) \ge 1$. A sufficient condition for (10) can be deduced as follows. Recursive substitution gives

$$\sigma_t^2 \le K + K \sum_{l=1}^{\infty} \left( \sum_{j_1=1}^{\infty} \cdots \sum_{j_l=1}^{\infty} \psi_{0j_1} \cdots \psi_{0j_l} \varepsilon_{t-j_1}^2 \varepsilon_{t-j_1-j_2}^2 \cdots \varepsilon_{t-j_1-\cdots-j_l}^2 \right),$$

so by the $c_r$-inequality,

$$\sigma_t^{2\rho} \le K + K \sum_{l=1}^{\infty} \left( \sum_{j_1=1}^{\infty} \cdots \sum_{j_l=1}^{\infty} \psi_{0j_1}^{\rho} \cdots \psi_{0j_l}^{\rho} |\varepsilon_{t-j_1}|^{2\rho} \right.$$
$$\left. \times |\varepsilon_{t-j_1-j_2}|^{2\rho} \cdots |\varepsilon_{t-j_1-\cdots-j_l}|^{2\rho} \right).$$

Thus, from Lemma 2,

$$E|x_t|^{2\rho} < E|\sigma_t|^{2\rho} \le K + K \sum_{l=0}^{\infty} \left( E|\varepsilon_0|^{2\rho} \sum_{j=1}^{\infty} \psi_{0j}^{\rho} \right)^l.$$

The last bound is finite if and only if

$$
(32) \qquad E|\varepsilon_0|^{2\rho} \sum_{j=1}^{\infty} \psi_{0j}^{\rho} < 1.
$$

Thus, (10) holds if there is a $\rho$ satisfying (9) and (32). Recursive substitution and the $c_r$-inequality were also used by Nelson ([27], Corollary) to upper-bound $E|\sigma_t|^{2\rho}$ in the GARCH(1, 1) case, but he employed the simple dynamic structure available there, and (32) does not reduce to his necessary and sufficient condition.

If $\psi(1; \zeta_0) < 1$, (32) adds nothing because we already know that $Ex_0^2 < \infty$ here, but if $\psi(1; \zeta_0) \geq 1$, the second factor on the left-hand side of (32) exceeds 1 and increases with $\rho$; the question is whether the first factor, which is less than 1 and decreases with $\rho$ [due to Assumption A($q$)], can over-compensate. Analytic verification of (32) for given $\zeta_0, \rho$ seems in general infeasible, and numerical verification highly problematic when the $\psi_j$ decay slowly. However, consider the family of densities

$$
(33) \qquad f(\varepsilon) = \exp[-\{\alpha(\gamma)|\varepsilon|\}^{1/\gamma}] / \{2\gamma\Gamma(\gamma)\alpha(\gamma)\}
$$

for $\gamma > 0$, where $\alpha(\gamma) = \{\Gamma(\gamma)/\Gamma(3\gamma)\}^{1/2}$ (also used by Nelson [28] to model the innovation of the exponential GARCH model). We have $E\varepsilon_0 = 0$, $E\varepsilon_0^2 = 1$ as necessary, Assumption A($q$) is satisfied for all $q > 0$, and $E|\varepsilon_0|^{2\rho} = \Gamma((2\rho + 1)\gamma)/\{\Gamma(\gamma)^{1-\rho}\Gamma(3\gamma)^{\rho}\}$. In case $\gamma = 0.5$, (33) is the normal density, for which $\hat{\theta}_T$ is asymptotically efficient. Here $E|\varepsilon_0|^{2\rho} = 2^{\rho}\Gamma(\rho + 0.5)/\sqrt{\pi}$, and numerical calculations for FIGARCH(0, $d_0$, 0) cast doubt on (32). In case $\gamma = 1$, (33) is the Laplace density, with $E|\varepsilon_0|^{2\rho} = 2^{\rho-1}\Gamma(2\rho + 1)$. As $\gamma$ increases, $E|\varepsilon_0|^{2\rho}$ can be made small for fixed $\rho < 1$, for example, with $\rho = 0.95$, it is 0.64 when $\gamma = 10$ and 0.42 when $\gamma = 20$.

**4. Generalized exponential and hyperbolic models.** FGARCH($n, d_0, m$) [and FIGARCH($n, d_0, m$)] processes require $d_0 \in (0, 1)$. For $d = 1$, (29) reduces to (23), and for $d > 1$, at least one coefficient in the expansion of (23) is negative, leading to the possibility of negative $\psi_j(\zeta)$. Because FGARCH $\psi_j(\zeta)$ decay like $j^{-d-1}$, a large mathematical gap is left relative to GARCH processes. Even if exponential decay is anticipated, there is a case for more direct modeling of the $\psi_j(\zeta)$ than provided by GARCH($n, m$), since it is the $\psi_j(\zeta)$ and their derivatives that must be formed in point and interval estimation based on the PMLE.

Consider the choices

$$
(34) \qquad \psi_j(\zeta) = \sum_{i=1}^{m} \Gamma(f_i + 1)^{-1} e_i d^{f_i+1} j^{f_i} e^{-dj},
$$

$$
(35) \qquad \psi_j(\zeta) = \sum_{i=1}^{m} \Gamma(f_i + 1)^{-1} e_i d \ln^{f_i}(j + 1)(j + 1)^{-d-1},
$$

where $d = d(\zeta)$ and the $e_i = e_i(\zeta)$, $f_i = f_i(\zeta)$ are such that $\Upsilon$ satisfies

$$(36) \qquad\qquad\qquad d \in (0, \infty),$$

$$(37) \qquad\qquad\qquad e_i > 0, \qquad i = 1, \ldots, m,$$

$$(38) \qquad\qquad\qquad 0 \le f_1 \le \cdots \le f_m < \infty,$$

with $2m + 1 \ge r$. Given (1)–(4) and (22), we call $x_t$ generated by (34) a generalized exponential, $\mathrm{GEXP}(m)$, process, and $x_t$ generated by (35) a generalized hyperbolic, $\mathrm{GHYP}(m)$, process. Condition (38) is sufficient but not necessary for $\psi_j(\zeta) > 0$, all $j \ge 1$. By choosing $m$ large enough in (34) or (35), any finite $\psi(1; \zeta)$ can be arbitrarily well approximated, but (34) and (35) can also achieve parsimony. For real $x \ge 1$, $x^f e^{-dx}$ and $(\ln x)^f x^{-d-1}$ decay monotonically if $f = 0$, and for $f > 0$, have single maxima at $f/d$ and $e^{f/(d+1)}$, respectively. Thus, with $m = 1$ and $f_1 = 0$, we have monotonic decay in (34) and (35); otherwise, both can exhibit lack of monotonicity, while eventually decaying exponentially or hyperbolically. The scale factors in (34) and (35) are so expressed because $x^f e^{-dx}$ and $(\ln x)^f x^{-d-1}$ integrate over $(0, \infty)$ to $\Gamma(f+1)/d^{f+1}$ and $\Gamma(f+1)/d$, respectively, so that $\psi(1; \zeta) \simeq \sum_{i=1}^{m} e_i$ in both cases, but the approximation may not be very close and the "integrated" case is less easy to distinguish than in GARCH and FGARCH models (though it would be possible to alternatively scale the weights by infinite sums to achieve equality).

The following corollary covers (34) and (35) simultaneously, and implies the special case when the $f_i$ are specified a priori, for example, to be non-negative integers; strictly speaking, when the true value of $f_1$ is unknown, Assumption C prevents it from being zero.

COROLLARY 2. *Let $\psi(z; \zeta)$ be given by (22) and (34) or (35) with $m \ge 1$ and let $d$ and the $e_i, f_i$ be continuously differentiable. For some $\delta > 0$, let Assumptions A$(2+\delta)$, B, C and E hold, with all $\zeta \in \Upsilon$ satisfying (36)–(38) and*

$$\mathrm{rank}\left\{ \frac{\partial}{\partial \zeta}(e_1, f_1, \ldots, e_m, f_m, d) \right\} = r.$$

*Then (15) is true. Let also $d$ and the $e_i, f_i$ be thrice continuously differentiable and Assumption A(4) hold, and $d_0 = d(\zeta_0) > \frac{1}{2}$ in case of (35). Then (17) is true.*

PROOF. Given (36)–(38) and the proofs of Corollaries 1 and 2, the verification of Assumptions D and F$(l)$ is straightforward. We check Assumption G for (35) only, a very similar type of proof holding for (34). We have

$$\psi_j^{(1)} = \begin{bmatrix} E(u'_{1j}, \ldots, u'_{mj})' \\ v_j \end{bmatrix} d(j+1)^{-d-1},$$

where

$$u_{ij} = (\ln\ln(j+1) - (\partial/\partial f_i)\ln\Gamma(f_{i+1}),1)'\ln^{f_i}(j+1), \qquad i = 1,\ldots,r,$$

$$v_j = -\sum_{i=1}^{m} e_i\Gamma(f_{i+1})^{-1}\ln^{f_{i+1}}(j+1),$$

and $E$ is the diagonal matrix whose $(2i-1)$st diagonal element is $e_i$, and whose even diagonal elements are all 1. Fixing $\zeta$, we show first that the leading $(r-1)\times(r-1)$ submatrix of $\Psi_{(j_1,\ldots,j_r)}(\zeta)$ has full rank, equivalently, that $U_m$ has full rank, where, for $i = 1,\ldots,m$, the $(2i)\times(2i)$ matrix $U_i$ has $(k,\ell)$th $2\times1$ sub-vector $u_{kj_\ell}$, $k = 1,\ldots,i$, $\ell = 1,\ldots,2i$. Suppose, for some $i = 1,\ldots,m-1$ and given $j_1,\ldots,j_{2i}$, that $U_i$ has full rank, and partition the rows and columns of $U_{i+1}$ in the ratio $2i:2$, calling its $(k,\ell)$th submatrix $U_{k\ell}$ (so $U_{11} = U_i$). Take $j_{2i+2} = j_{2i+1}^2$. Because $\ln\ln x$ strictly increases in $x > 1$, it follows that $U_{22}$ is nonsingular and $\|U_{22}^{-1}\| = O(\ln\ln j_{2i+1}\ln^{-f_{i+1}}j_{2i+1})$. Noting that $\|U_{12}\| = O(\ln\ln j_{2i+1}\ln^{f_i}j_{2i+1})$, while $U_{11}$ and $U_{21}$ depend only on $j_1,\ldots,j_{2i}$, we can choose $j_{2i+1}$ such that $U_{11} - U_{12}U_{22}^{-1}U_{21}$ differs negligibly from $U_{11}$. Thus, $U_{i+1}$ has full rank. Since, for $f_1 \geq 0$, $U_1$ has full rank (e.g., when $j_1 = 1$, $j_2 = 2$), it follows by induction that $U_m$ has full rank. Since $v_j$ is dominated by a term of order $\ln^{f_m+1}j$, while $\|u_{ij}\| = O(\ln\ln j\ln^{f_i}j)$, a similar argument shows that $j_r$ can then be chosen large enough, to complete verification of Assumption G.   $\square$

## 5. Technical lemmas. Define

$$\sigma_t^{*2}(\theta) = \omega + \sum_{j=1}^{\infty}\psi_j(\zeta)x_{t-j}^2, \qquad \sigma_t^{*2} = \omega_U + \sum_{j=1}^{\infty}\sup_{\zeta\in\Upsilon}\psi_j(\zeta)x_{t-j}^2.$$

LEMMA 1.   *Under Assumptions* B *and* D, *for all* $\theta\in\Theta$, $t\in\mathbb{Z}$,

$$K^{-1}\sigma_t^{*2}(\theta) \leq \sigma_t^2(\theta) \leq K\sigma_t^{*2}(\theta) \qquad a.s.$$

PROOF.   A simple extension of [21], Lemma 1.   $\square$

LEMMA 2.   *Under Assumptions* A(2), B, C, D *and* E, *for all* $t\in\mathbb{Z}$,

$$\tag{39} E|x_t|^{2\rho} < E\sigma_t^{2\rho} \leq E\sup_{\theta\in\Theta}\sigma_t^{2\rho}(\theta) \leq KE\sigma_t^{*2\rho} \leq KE|x_t|^{2\rho} \leq K,$$

$$\tag{40} \inf_{\theta\in\Theta}\sigma_t^2(\theta) > 0, \qquad \sup_{\theta\in\Theta}\sigma_t^2(\theta) < K\sigma_t^{*2} < \infty \qquad a.s.,$$

$$\tag{41} E\sup_{\theta\in\Theta}|\ln\sigma_t^2(\theta)| \leq K.$$

PROOF.  With respect to (39), the first inequality follows from Jensen's inequality, the second is obvious, the third follows from Lemma 1, the fourth follows from the $c_r$-inequality, (7) and (9), while the last one is due to (10). The proof of (40) uses Lemma 1, $\sigma_t^2(\theta) \geq \omega_L$, (10) and [23], page 121. To prove (41), $|\ln x| \leq x + x^{-1}$ for $x > 0$ and Lemma 2 give

$$E \sup_{\theta \in \Theta} |\ln \sigma_t^2(\theta)| \leq \rho^{-1} E \sup_{\theta \in \Theta} \sigma_t^{2\rho}(\theta) + E \left\{ \inf_{\theta \in \Theta} \sigma_t^2(\theta) \right\}^{-1} \leq K.$$

□

LEMMA 3.  *Under Assumptions* D, E *and* F(l), *for all* $\theta \in \Theta$, $\sigma_t^2(\theta)$, $q_t(\theta)$ *and their first* $l$ *derivatives are strictly stationary and ergodic.*

PROOF.  Follows straightforwardly from the assumptions.  □

LEMMA 4.  *Under Assumption* A(2), *for positive integer* $k < (b+1)n/2$,

$$E \left( \sum_{t=1}^n \varepsilon_t^2 \right)^{-k} < \infty. \tag{42}$$

PROOF.  Denote by $M_X(t) = E(e^{tX})$ the moment-generating function of a random variable $X$. By Cressie et al. [8], the left-hand side of (42) is proportional to

$$\int_0^\infty t^{k-1} M_{\sum \varepsilon_t^2}(-t) \, dt = \int_0^\infty t^{k-1} M_{\varepsilon_0^2}^n(-t) \, dt$$

$$\leq \int_0^1 t^{k-1} \, dt + \int_1^\infty t^{k-1} M_{\varepsilon_0^2}^n(-t) \, dt. \tag{43}$$

It suffices to show that the last integral is bounded. For all $\delta > 0$, there exists $\eta > 0$ such that $L(\varepsilon^{-1}) \leq \varepsilon^{-\delta}$, $\varepsilon \in (0, \eta)$, so

$$M_{\varepsilon_0^2}(-t) = \int_{-\infty}^\infty e^{-t\varepsilon^2} f(\varepsilon) \, d\varepsilon \leq K \int_0^\eta e^{-t\varepsilon^2} \varepsilon^{b-\delta} \, d\varepsilon + 2e^{-t\eta^2}.$$

The last integral is bounded by

$$K t^{(\delta-b-1)/2} \int_0^\infty e^{-\varepsilon} \varepsilon^{(\delta-b-1)/2} \, d\varepsilon \leq K t^{(\delta-b-1)/2}.$$

Thus, (43) is finite if $k + n(\delta - b - 1)/2 < 0$, that is, since $\delta$ is arbitrary, if $k < (b+1)n/2$.  □

The previous version of the paper included a longer, independently obtained, proof of the following lemma which we have been able to shorten in one respect by using an idea of Berkes, Horváth and Kokoszka [3] in a corresponding lemma covering the GARCH$(n, m)$ case.

LEMMA 5. *Under Assumptions* A$(q)$, B, C *and* D, *for* $p < q/2$,

$$E \sup_{\theta \in \Theta} \left( \frac{\sigma_t^2}{\sigma_t^2(\theta)} \right)^p \leq K < \infty.$$

PROOF. We have

$$\sigma_t^2 = \omega_0 + \psi_{01} x_{t-1}^2 + \sum_{j=2}^{\infty} \psi_{0j} x_{t-j}^2 \leq \omega_0 + \psi_{01} \sigma_{t-1}^2 \varepsilon_{t-1}^2 + K \sigma_{t-1}^2$$

from (8). Thus, $\sigma_t^2 / \sigma_{t-1}^2 \leq K(1 + \varepsilon_{t-1}^2)$ and thence, for fixed $j \geq 1$, $\sigma_t^2 / \sigma_{t-j}^2 \leq K h_{tj}$, where $h_{tj} = \prod_{i=1}^{j}(1 + \varepsilon_{t-i}^2)$. For any $M < \infty$,

$$\frac{\sigma_t^2}{\sigma_t^2(\theta)} \leq \frac{K \sigma_t^2}{\sigma_t^{*2}(\theta)} \leq K \left( \frac{\omega}{\sigma_t^2} + \sum_{j=1}^{M} \psi_j(\zeta) \varepsilon_{t-j}^2 \frac{\sigma_{t-j}^2}{\sigma_t^2} \right)^{-1}$$

$$\leq \frac{K h_{tM} / \{ \inf_{\zeta \in \Upsilon} \inf_{j=1,\dots,M} \psi_j(\zeta) \}}{\sum_{j=1}^{M} \varepsilon_{t-j}^2}.$$

The proof can now be completed much as in the proof of Lemma 5.1 of [3], using Hölder's inequality as there but employing our Lemma 4 and taking $M > 2pq/[(b+1)(q-2p)]$. □

LEMMA 6. *Under Assumptions* A$(2)$, B, C, D, E *and* F$(l)$, *for all* $p > 0$ *and* $k \leq l$,

$$E \sup_{\theta \in \Theta} \left| \frac{1}{\sigma_t^2(\theta)} \frac{\partial^k \sigma_t^2(\theta)}{\partial \theta_{i_1} \cdots \partial \theta_{i_k}} \right|^p < \infty, \tag{44}$$

$$E \sup_{\theta \in \Theta} \left| \frac{1}{\bar{\sigma}_t^2(\theta)} \frac{\partial^k \bar{\sigma}_t^2(\theta)}{\partial \theta_{i_1} \cdots \partial \theta_{i_k}} \right|^p < \infty. \tag{45}$$

PROOF. Take $i_1 \leq i_2 \leq \cdots \leq i_k$. First assume $i_1 \geq 3$, whence, for given $k$ and $i_1, \dots, i_k$,

$$\frac{\partial^k \sigma_t^2(\theta)}{\partial \theta_{i_1} \cdots \partial \theta_{i_k}} = \sum_{j=1}^{\infty} \xi_j(\zeta) x_{t-j}^2(\mu),$$

where $\xi_j(\zeta) = \partial^k \psi_j(\zeta) / \partial \zeta_{i_1-2} \cdots \partial \zeta_{i_k-2}$. Now

$$\left| \sum_{j=1}^{\infty} \xi_j(\zeta) x_{t-j}^2(\mu) \right| \leq 2 \sum_{j=1}^{\infty} |\xi_j(\zeta)|(x_{t-j}^2 + \mu^2),$$

so using Lemma 1,

$$\left| \frac{1}{\sigma_t^2(\theta)} \frac{\partial^k \sigma_t^2(\theta)}{\partial \theta_{i_1} \cdots \partial \theta_{i_k}} \right| \leq \frac{2 \sum_{j=1}^{\infty} |\xi_j(\zeta)| x_{t-j}^2}{\sigma_t^{*2}(\theta)} + K \sum_{j=1}^{\infty} |\xi_j(\zeta)|.$$

It suffices to take $p > 1$. By Hölder's inequality,

$$\sum_{j=1}^{\infty} |\xi_j(\zeta)||x_{t-j}^2 \le \left\{\sum_{j=1}^{\infty} |\xi_j(\zeta)|^{p/\rho} \psi_j(\zeta)^{1-p/\rho} x_{t-j}^2 \right\}^{\rho/p} \left\{\sum_{j=1}^{\infty} \psi_j(\zeta) x_{t-j}^2 \right\}^{1-\rho/p},$$

so

$$\left\{\frac{\sum_{j=1}^{\infty} |\xi_j(\zeta)||x_{t-j}^2}{\sigma_t^{*2}(\theta)}\right\}^p \le K \sum_{j=1}^{\infty} |\xi_j(\zeta)|^p \psi_j(\zeta)^{\rho-p} |x_{t-j}|^{2\rho}.$$

By Assumption F($l$), for all $\eta > 0$

$$\sup_{\zeta \in \Upsilon} |\xi_j(\zeta)|^p \psi_j(\zeta)^{\rho-p} \le K \sup_{\zeta \in \Upsilon} \psi_j(\zeta)^{\rho-\eta p} \le K j^{-(\underline{d}+1)(\rho-\eta p)}.$$

Since $\rho(\underline{d}+1) > 1$, we may choose $\eta$ such that $(\underline{d}+1)(\rho-p\eta) > 1$. Thus,

$$E \sup_{\theta \in \Theta} \left\{\frac{\sum_{j=1}^{\infty} |\xi_j(\zeta)||x_{t-j}^2}{\sigma_t^{*2}(\theta)}\right\}^p < \infty.$$

The above proof implies that also

$$\sup_{\zeta \in \Upsilon} \left\{\sum_{j=1}^{\infty} |\xi_j(\zeta)| \right\}^p < \infty,$$

whence, the proof of (44) with $i_1 \ge 3$ is concluded. Next take $i_1 = 2$. If $i_2 > 2$,

$$(46) \qquad \frac{\partial^k \sigma_t^2(\theta)}{\partial \theta_{i_1} \cdots \partial \theta_{i_k}} = -2 \sum_{j=1}^{\infty} \xi_j(\zeta) x_{t-j}(\mu),$$

where now $\xi_j(\zeta) = \partial^{k-1} \psi_j(\zeta)/\partial \zeta_{i_2-2} \cdots \partial \zeta_{i_k-2}$, while if $i_2 = 2$, $i_3 > 2$,

$$\frac{\partial^k \sigma_t^2(\theta)}{\partial \theta_{i_1} \cdots \partial \theta_{i_k}} = -2 \sum_{j=1}^{\infty} \xi_j(\zeta),$$

where now $\xi_j(\zeta) = \partial^{k-2} \psi_j(\zeta)/\partial \zeta_{i_3-2} \cdots \partial \zeta_{i_k-2}$. In the first of these cases the proof is seen to be very similar to that above after noting that, by the Cauchy inequality, (46) is bounded by

$$K \left\{\sum_{j=1}^{\infty} |\xi_j(\zeta)| x_{t-j}^2 \sum_{j=1}^{\infty} |\xi_j(\zeta)| \right\}^{1/2} + K \sum_{j=1}^{\infty} |\xi_j(\zeta)|,$$

while in the second it is more immediate; we thus omit the details. We are left with the cases $i_1 = i_2 = i_3 = 2$ and $i_1 = 1$, both of which are trivial. The details for (45) are very similar (the truncations in numerator and denominator match), and are thus omitted.  $\square$

Define

$$g_t(\theta) = u_t(\theta)u_t'(\theta), \qquad G_T(\theta) = T^{-1}\sum_{t=1}^{T} g_t(\theta).$$

LEMMA 7. *For some $\delta > 0$, under Assumptions* A$(2+\delta)$, B, C, D, E *and* F*(1)*,

$$(47) \qquad \sup_{\theta \in \Theta}|Q_T(\theta) - Q(\theta)| \to 0 \qquad a.s. \ as \ T \to \infty,$$

*and $Q(\theta)$ is continuous in $\theta$. If also Assumption* F*(2) holds,*

$$(48) \qquad \sup_{\theta \in \Theta}\|G_T(\theta) - G(\theta)\| \to 0 \qquad a.s. \ as \ T \to \infty,$$

*and $G(\theta)$ is continuous in $\theta$. If also Assumption* F*(3) holds,*

$$(49) \qquad \sup_{\theta \in \Theta}\|H_T(\theta) - H(\theta)\| \to 0 \qquad a.s. \ as \ T \to \infty,$$

*and $H(\theta)$ is continuous in $\theta$.*

PROOF. To prove (47), note first that, by Lemmas 1, 2, 3 and 5,

$$\sup_{\Theta} E|q_0(\theta)| \le \sup_{\Theta} E|\log \sigma_0^2(\theta)| + \sup_{\Theta} E\chi_0(\theta) < \infty.$$

Thus, by ergodicity

$$Q_T(\theta) \to Q(\theta) \qquad \text{a.s.},$$

for all $\theta \in \Theta$. Then uniform convergence follows on establishing the equicontinuity property

$$\sup_{\tilde{\theta}\,:\,\|\tilde{\theta}-\theta\|<\varepsilon}|Q_T(\tilde{\theta}) - Q_T(\theta)| \to 0 \qquad \text{a.s.},$$

as $\varepsilon \to 0$, and continuity of $Q(\theta)$. By the mean value theorem it suffices to show that

$$\sup_{\Theta}\left\|\frac{\partial Q_T(\theta)}{\partial \theta}\right\| + \sup_{\Theta}\left\|\frac{\partial Q(\theta)}{\partial \theta}\right\| < \infty \qquad \text{a.s.},$$

which, by Loève ([23], page 121) and identity of distribution, is implied by $E\sup_{\Theta}\|u_0(\theta)\| < \infty$. By the definition of $u_t(\theta)$, and $x_t^2(\mu) \le K(x_t^2+1)$, $\|\nu_t(\theta)\| \le 2(|x_t|+1)$, we have

$$\|u_t(\theta)\| \le K\left[\|\tau_t(\theta)\|\left\{1 + \varepsilon_t^2\frac{\sigma_t^2}{\sigma_t^2(\theta)}\right\} + |\varepsilon_t|\frac{\sigma_t}{\sigma_t(\theta)} + 1\right].$$

Thus, $E \sup_{\Theta} \|u_0(\theta)\|$ is bounded by a constant times

$$E \sup_{\Theta} \|\tau_0(\theta)\| + \left[ E \sup_{\Theta} \left\{ \frac{\sigma_0^2}{\sigma_0^2(\theta)} \right\}^p \right]^{1/p} \left[ E \sup_{\Theta} \|\tau_0(\theta)\|^{p/(p-1)} \right]^{1-1/p}$$

$$+ E \sup_{\Theta} \left\{ \frac{\sigma_0}{\sigma_0(\theta)} \right\} + 1$$

for all $p > 1$. On choosing $p < 1 + \delta/2$, this is finite by Lemmas 5 and 6. (Our use of Lemmas 5 and 6 is similar to Berkes, Horváth and Kokoszka's [3] use of their Lemmas 5.1 and 5.2 in the GARCH($n, m$) case.) This completes the proof of (47). Then (48) and (49) follow by applying analogous arguments to those above, and so we omit the details; indeed, (48) and (49) are only used in the proof of consistency of $\bar{G}_T(\hat{\theta}_T)$, $\bar{H}_T(\hat{\theta}_T)$ for $G_0, H_0$, where convergence over only a neighborhood of $\theta_0$ would suffice. □

Lemma 8. *Under Assumptions* A$(2 + \delta)$, B, C, D, E *and* F$(1)$,

$$\tag{50} \sup_{\theta \in \Theta} |Q_T(\theta) - \bar{Q}_T(\theta)| \to 0 \qquad a.s. \ as \ T \to \infty.$$

*If also Assumption* F$(2)$ *holds,*

$$\tag{51} \sup_{\theta \in \Theta} \|G_T(\theta) - \bar{G}_T(\theta)\| \to 0 \qquad a.s. \ as \ T \to \infty.$$

*If also Assumption* F$(3)$ *holds,*

$$\tag{52} \sup_{\theta \in \Theta} \|H_T(\theta) - \bar{H}_T(\theta)\| \to 0 \qquad a.s. \ as \ T \to \infty.$$

Proof. We have $\hat{Q}_T(\theta) - Q_T(\theta) = A_T(\theta) + B_T(\theta)$, where

$$A_T(\theta) = T^{-1} \sum_{t=1}^{T} \ln \left[ \frac{\bar{\sigma}_t^2(\theta)}{\sigma_t^2(\theta)} \right], \qquad B_T(\theta) = T^{-1} \sum_{t=1}^{T} x_t^2(\mu) \{ \bar{\sigma}_t^{-2}(\theta) - \sigma_t^{-2}(\theta) \}.$$

Because

$$\sigma_t^2(\theta) = \bar{\sigma}_t^2(\theta) + \sum_{j=0}^{\infty} \psi_{j+t}(\zeta) x_{-j}^2(\mu),$$

$\ln(1 + x) \leq x$ for $x > 0$ and $\sigma_t^2(\theta) \geq \omega_L > 0$, it follows that

$$|A_T(\theta)| \leq K T^{-1} \sum_{t=1}^{T} \{ \sigma_t^2(\theta) - \bar{\sigma}_t^2(\theta) \}$$

$$\leq K T^{-1} \sum_{t=1}^{T} \sum_{j=t}^{\infty} \psi_j(\zeta) x_{t-j}^2(\mu)$$

$$\leq K T^{-1} \sum_{t=0}^{\infty} \left\{ \sum_{j=t+1}^{t+T} \psi_j(\zeta) \right\} x_{-t}^2(\mu).$$

Now from (7),

$$\sup_{\zeta \in \Upsilon} \sum_{j=t+1}^{t+T} \psi_j(\zeta) \le K \sum_{j=t+1}^{t+T} j^{-\underline{d}-1} \le K \min(t+1, T)(t+1)^{-\underline{d}-1}.$$

Thus,

$$\text{(53)} \qquad \sup_{\Theta} A_T(\theta) \le K T^{-1} \sum_{t=0}^{T} (t+1)^{-\underline{d}} (x_{-t}^2 + 1) + K \sum_{t=T}^{\infty} t^{-\underline{d}-1}(x_{-t}^2 + 1).$$

From the $c_r$-inequality, (9) and (10), $\sum_{t=1}^{\infty}(t+1)^{-\underline{d}-1} x_{-t}^2$ has finite $\rho$th moment, and thus, by Loève ([23], page 121), is a.s. finite. Thus, the second term of (53) tends to zero a.s. as $T \to \infty$, while the first does so for the same reasons combined with the Kronecker lemma. Next,

$$\text{(54)} \qquad \begin{aligned} |B_T(\theta)| &\le K T^{-1} \sum_{t=1}^{T} \chi_t(\theta) \sum_{j=t}^{\infty} \psi_j(\zeta) x_{t-j}^2(\mu) \\ &\le K T^{-1} \sum_{t=1}^{T} \chi_t(\theta) \sum_{j=t}^{\infty} j^{-\underline{d}-1}(x_{t-j}^2 + 1). \end{aligned}$$

From previous remarks, $\sum_{j=t}^{\infty} j^{-\underline{d}-1}(x_{t-j}^2 + 1) \to 0$ a.s. Also, for each $\theta$, a.s.

$$T^{-1} \sum_{t=1}^{T} \chi_t(\theta) \to E\chi_0(\theta) \le K \left\{ E\left( \frac{\sigma_0^2}{\sigma_0^2(\theta)} \right) + 1 \right\} \le K$$

by ergodicity and Lemma 5. Thus, (54) $\to 0$ a.s. by the Toeplitz lemma. The convergence is uniform in $\theta$ because, from the proof of Lemma 7, for all $\theta \in \Theta$,

$$\sup_{\tilde\theta : \|\tilde\theta - \theta\| < \varepsilon} \|\chi_0(\tilde\theta) - \chi_0(\theta)\| \to 0 \qquad \text{a.s.,}$$

as $\varepsilon \to 0$. This completes the proof of (50). We omit the proofs of (51) and (52) as they involve the same kind of arguments. $\square$

LEMMA 9. *For some $\delta > 0$, under Assumptions* A$(2 + \delta)$, B, C, D, E, F$(1)$ *and* G, *$M(\theta)$ is finite and positive definite for all $\theta \in \Theta$.*

PROOF. Fix $\theta \in \Theta$. Finiteness of $M(\theta)$ follows from Lemma 6. Positive definiteness follows (by an argument similar to that of Lumsdaine [24] in the GARCH$(1,1)$ case) if, for all nonnull $(r+2) \times 1$ vectors $\lambda$, $\lambda' M(\theta)\lambda = E\{\lambda' \tau_0(\theta)\}^2 > 0$, that is, that

$$\text{(55)} \qquad \lambda' \tau_0(\theta)\sigma_0^2(\theta) \neq 0 \qquad \text{a.s.,}$$

since $0 < \sigma_0^2(\theta) < \infty$ a.s. Define

$$\tau_{t\omega}(\theta) = \frac{\partial}{\partial\omega} \ln \sigma_t^2(\theta) = \sigma_t^{-2}(\theta),$$

$$\tau_{t\mu}(\theta) = \frac{\partial}{\partial\gamma} \ln \sigma_t^2(\theta) = -2\sigma_t^{-2}(\theta) \sum_{j=1}^{\infty} \psi_j(\zeta) x_{t-j}(\mu),$$

$$\tau_{t\zeta}(\theta) = \frac{\partial}{\partial\zeta} \ln \sigma_t^2(\theta) = \sigma_t^{-2}(\theta) \sum_{j=1}^{\infty} \psi_j^{(1)}(\zeta) x_{t-j}^2(\mu),$$

so that $\tau_t(\theta) = (\tau_{t\omega}(\theta),\ \tau_{t\mu}(\theta),\ \tau_{t\zeta}'(\theta))'$. Write $\lambda = (\lambda_1, \lambda_2, \lambda_3')'$, where $\lambda_1$ and $\lambda_2$ are scalar and $\lambda_3$ is $r \times 1$. Consider first the case $\lambda_1 = \lambda_2 = 0$, $\lambda_3 \neq 0$. Suppose (55) does not hold. Then we must have

$$\sum_{j=1}^{\infty} \lambda_3' \psi_j^{(1)}(\zeta) x_{t-j}^2(\mu) = 0 \qquad \text{a.s.}$$

If $\lambda_3' \psi_1^{(1)}(\zeta) \neq 0$, it follows that

$$(56) \qquad (\sigma_{t-1}\varepsilon_{t-1} + \mu_0 - \mu)^2 = -\{\lambda_3'\psi_j^{(1)}(\zeta)\}^{-1} \sum_{j=2}^{\infty} \lambda_3'\psi_j^{(1)}(\zeta) x_{-j}^2(\mu).$$

Since $\sigma_{t-1} > 0$ a.s., the left-hand side involves the nondegenerate random variable $\varepsilon_{t-1}$, which is independent of the right-hand side, so (56) cannot hold. Thus, $\lambda_3' \psi_1^{(1)}(\zeta) = 0$. Repeated application of this argument indicates that, for all $\zeta$, $\lambda_3' \psi_j^{(1)}(\zeta) = 0$, $j = 1, \ldots, j_r(\zeta)$. This is contradicted by Assumption G, so (56) cannot hold. Next consider the case $\lambda_1 = 0$, $\lambda_2 \neq 0$, $\lambda_3 = 0$. If (56) does not hold, we must have

$$(57) \qquad \sum_{j=1}^{\infty} \psi_j(\zeta) x_{t-j}(\mu) = 0 \qquad \text{a.s.}$$

Let $k$ be the smallest integer such that $\psi_k(\zeta) \neq 0$. Then (57) implies

$$\varepsilon_{t-k} = \sigma_{t-k}^{-1}(\theta) \left\{ \mu - \mu_0 - \psi_k^{-1}(\zeta) \sum_{j=k+1}^{\infty} \psi_j(\theta) x_{t-j}(\mu) \right\}.$$

But the left-hand side is nondegenerate and independent of the right-hand side, so (57) cannot hold. Next consider the case $\lambda_1 = 0$, $\lambda_2 \neq 0$, $\lambda_3 \neq 0$. If (55) is not true, then, taking $\lambda_2 = 1$, we must have

$$(58) \qquad \sum_{j=1}^{\infty} \{\lambda_3'\psi_j^{(1)}(\zeta) x_{t-j}(\mu) - 2\psi_j(\zeta)\} x_{t-j}(\mu) = 0 \qquad \text{a.s.}$$

Let $k$ be the smallest integer such that either $\lambda_3' \psi_k^{(1)}(\zeta) \neq 0$ or $\psi_k(\zeta) \neq 0$; the preceding argument indicates that there exists such $k$. Then we have

$$\{2\psi_k(\zeta) - \lambda_3' \psi_k^{(1)}(\zeta)(\sigma_{t-k}\varepsilon_{t-k} + \mu_0 - \mu)\}\{\sigma_{t-k}\varepsilon_{t-k} + \mu_0 - \mu\}$$
$$= \sum_{j=k+1}^{\infty} \{\lambda_3' \psi_j^{(1)}(\zeta) x_{t-j}(\mu) - 2\psi_j(\zeta)\} x_{t-j}(\mu) \qquad \text{a.s.}$$

The left-hand side is a.s. nonzero and involves the nondegenerate random variable $\varepsilon_{t-k}$, which is independent of the right-hand side, so (58) cannot hold. We are left with the cases where $\lambda_1 \neq 0$. Taking $\lambda_1 = -1$ and noting that $\sigma_t^2(\theta)\tau_{t\omega}(\theta) \equiv 1$, the preceding arguments indicate that there exist no $\lambda_2$ and $\lambda_3$ such that

$$\lambda_2 \sigma_t^2(\theta)\tau_{t\mu}(\theta) + \lambda_3' \sigma_t^2(\theta)\tau_{t\zeta}(\theta) = 1 \qquad \text{a.s.}$$

$\square$

LEMMA 10. *For some $\delta > 0$, under Assumptions* A$(2+\delta)$, B, C, D, E, F$(1)$ *and* H,

$$\inf_{\substack{\theta \in \Theta \\ \theta \neq \theta_0}} Q(\theta) > Q(\theta_0).$$

PROOF. We have

$$Q(\theta) - Q(\theta_0) = E\left[\frac{\sigma_0^2}{\sigma^2(\theta)} - \ln\left\{\frac{\sigma_0^2}{\sigma^2(\theta)}\right\} - 1\right] + (\mu - \mu_0)^2 E\left[\frac{1}{\sigma_0^2(\theta)}\right].$$

The second term on the right-hand side is zero only when $\mu = \mu_0$ and is positive otherwise. Because $x - \ln x - 1 \geq 0$ for $x > 0$, with equality only when $x = 1$, it remains to show that

$$(59) \qquad \ln \sigma_0^2(\theta) = \ln \sigma_0^2 \qquad \text{a.s., some } \theta \neq \theta_0.$$

By the mean value theorem, (59) implies that $(\theta - \theta_0)' \tau_0(\bar{\theta}) = 0$ a.s., for $\theta \neq \theta_0$ and some $\bar{\theta}$ such that $\|\bar{\theta} - \theta_0\| \leq \|\theta - \theta_0\|$. But by Lemma 9 there is no such $\bar{\theta}$. $\square$

**Acknowledgments.** We thank an Associate Editor and referees for a number of helpful comments that have led to a considerable improvement in the paper, and Fabrizio Iacone for help with the numerical calculations referred to in Section 3.

Department of Economics
London School of Economics
Houghton Street
London WC2A 2AE
United Kingdom
E-mail: p.m.robinson@lse.ac.uk

Tanaka Business School
Imperial College London
South Kensington Campus
London SW7 2AZ
United Kingdom
E-mail: p.zaffaroni@imperial.ac.uk